\documentclass[12pt]{amsart}
\usepackage{amsmath,amsfonts,amssymb,amsthm}
\usepackage{pinlabel}
\usepackage{graphicx}
\usepackage{mathtools}
\usepackage[all]{xy}
\usepackage{hyperref}
\usepackage[usenames,dvipsnames]{color}
\usepackage{color}
\xyoption{dvips}

\newtheorem{theorem}{Theorem}[section]

\newtheorem{corollary}[theorem]{Corollary}

\theoremstyle{remark}
\newtheorem{remark}[theorem]{Remark}
\theoremstyle{definition}

\newtheorem{ex}[theorem]{Example}

\newcommand{\R}{\mathbb{R}}

\numberwithin{equation}{section}

\begin{document}

\title[Non-regular Lagrangian concordances between fillable Legendrian knots]{Non-regular Lagrangian concordances between Lagrangian fillable Legendrian knots}

\author{Georgios Dimitroglou Rizell}
\author{Roman Golovko}

\begin{abstract}
In this short note, we construct a family of non-regular, and therefore non-decomposable, Lagrangian concordances between
Lagrangian fillable Legendrian knots in the standard contact $3$-dimensional sphere.  More precisely, for every decomposable Lagrangian concordance from $\Lambda_-$ to $\Lambda_+$, where $\Lambda_{\pm}$ are smoothly non-isotopic Legendrian knots, we construct a non-regular Lagrangian concordance from $\Lambda'_+$ to $\Lambda'_-$, where both $\Lambda'_\pm$ are Lagrangian fillable Legendrian knots obtained from $\Lambda_{\pm}$ by sufficiently many positive and negative stabilisations, followed by a sequence of Legendrian satellite operations.
\end{abstract}

\address{Department of Mathematics, Uppsala University,  Uppsala, Sweden}
\email{georgios.dimitroglou@math.uu.se}

\address{Faculty of Mathematics and Physics, Charles University, Prague, Czech Republic} \email{golovko@karlin.mff.cuni.cz}
\date{\today}
\subjclass[2010]{Primary 53D12; Secondary 53D42}

\keywords{Legendrian submanifold, decomposable Lagrangian concordance, regular Lagrangian concordance, Lagrangian filling}
\thanks{The first author is supported by the Knut and Alice Wallenberg Scholar grant KAW 2023.0294 and the Swedish Research Council Centre of Excellence grant VR 2022-06593: The Centre of Excellence in Geometry and Physics at Uppsala University. The second author is supported by the GA\v{C}R grant 26-20231L}

\maketitle

\section{Introduction and main results}

Weinstein cobordisms play a role in contact and symplectic topology analogous to the role that Morse cobordisms play in differential topology: they provide a structured way to build, modify, and relate contact and symplectic manifolds by attaching and modifying standard model handles that are compatible with a gradient-like flow.  More precisely,  Weinstein cobordism is a quadruple $(W,\omega,X,\Phi)$ such that $(W,\omega=d\beta)$ is an exact symplectic cobordism; $X$ is the Liouville vector field, which is defined to be a complete vector field given by $i_X \omega=\beta$ that enters the negative end of $W$ and exits from the positive end of $W$, and $\Phi$ is a Morse function on $W$ such that $X$ is a pseudo-gradient for $\Phi$. The existence of $\Phi$ allows us to get a decomposition of $(W,\omega,X,\Phi)$ into symplectic handles. For more details, we refer the reader to \cite{EliashbergWeinstein}. In stark contrast to the theory of ordinary smooth cobordisms, it is not clear when a Weinstein structure exists. There are definitely cases when the topology obstructs an existence; note that the Weinstein handles in dimension $2n$ all have index at most equal to $n$. One very entral open question in symplectic topology is understanding the mechanism that ensures the existence of a Weinstein structure.

The study of Lagrangians submanifolds occupies a prominent role in symplectic topology. Important such examples are Lagrangian cobordisms between Legendrian knots in Weinstein cobordisms, which are relevant objects both for the study of contact topology as well as for symplectic topology. Considerable progress has been made in constructing these cobordisms \cite{GeogrBourSabloffTraynor,LagrConcLegKnots,LegKnotsLagrCob,EtnyreLeversonribbon}, identifying obstructions to their existence \cite{GeogrBourSabloffTraynor,LagrConcLegKnots,FlThLagCob,ObstrLagrConc,LegKnotsLagrCob}, and distinguishing inequivalent examples \cite{CasalsGaoInfinite,CasalsNgBraidloops,CasalsCapovilla-SearleNewtonPolytopes,InfNumExLagrFil}. They also arise naturally in symplectic field theory \cite{IntroductionSFT}, as every Lagrangian cobordism induces a map between the contact homologies of its Legendrian ends. In particular, the most basic instance of such a cobordism are the orientable Lagrangian cobordisms of genus zero $L$ between Legendrian knots, i.e.  $L\hookrightarrow (W,\omega,X,\Phi)$ is an embedded, orientable Lagrangian genus zero surface, which is cylindrical over Legendrian knots at the positive and negative ends of $W$; such cobordisms are called Lagrangian concordances. It is necessary to note the fundamental difference between smooth and Lagrangian concordances. While smooth concordance is a symmetric relation on smooth knots, the Lagrangian concordance relation is directional and hence, as demonstrated by Chantraine \cite{LagrConcNotSymRel}, is not symmetric. For a background on Legendrian knots we refer the reader to \cite{EtnyreLegTransKnots, EtnyreNgSurvey}.

There is a particularly nice class of Lagrangians inside Weinstein manifolds whose topology can be reduced to understanding a handle-decomposition of its complement; these are the so-called regular Lagrangians, which \color{black} were introduced by Eliashberg-Lazarev-Murphy in  \cite{FlixibleLagrangians}. 
Recall that we say that a Lagrangian cobordism $L$ in  a Weinstein cobordism $(W,\omega,X,\Phi)$ 
is {\em regular} if $(W,\omega,X,\Phi)$  can be deformed to a
Weinstein structure $(W,\omega',X',\Phi')$ through Weinstein structures for which $L$ is an exact Lagrangian, and such that the new Liouville vector field $X'$ is tangent to $L$.  As for the question about existence of Weinstein structures, we currently also lack mechanisms for telling whether some given Lagrangian is regular or not. The questions whether e.g.~closed exact Lagrangian submanifolds of Weinstein domains all are regular is one of the most central open questions for understanding the symplectic topology of Lagrangians.\color{black}

Regular Lagrangian cobordisms are closely related to decomposable Lagrangian cobordisms as defined by Chantraine in \cite{Non-collarableSlices} and Ekholm--Honda--Kalman in \cite{LegKnotsLagrCob}. {\em Decomposable Lagrangian cobordisms} are those cobordisms that are obtained from a Lagrangian cobordism realization of a Legendrian isotopy, 
isolated standard unknot birth (induces Lagrangian 0-handle attachment), and Legendrian surgery (induces Lagrangian 1-handle attachment).  
From the work of Conway--Etnyre--Tosun  \cite[Theorem 1.10]{NonregImplNondecomp} it folllows that every decomposable Lagrangian cobordism is regular. 

The first examples of non-regular Lagrangian cobordisms
have been constructed by  Lin in dimension four \cite{LinExLagrCaps} and in high dimensions  by  Eliashberg-Murphy \cite{EliashbergMurphyLagCaps}. The examples of Lin and  of Eliashberg–Murphy have  empty positive ends  and stabilised (or loose) negative ends, in other words, they are Lagrangian caps of stabilised (or loose) Legendrians.  Recently, the first examples of non-regular Lagrangian concordances between Legendrian knots appeared in the work of the authors \cite{NonregStabilizedConc}. The key properties of the non-regular concordances appearing in \cite{NonregStabilizedConc} are that both positive and negative Legendrian ends of those concordances are  sufficiently stabilised and that the positive end is an unknot. That result was used in \cite{NonregStabilizedConc} to show that the relation of Lagrangian concordance between Legendrian knots is not a partial order; the same result for high-dimensional Legendrians has been proven by the second author in \cite{Notpartordhighdim}.
In this short note we extend the class of known non-regular Lagrangian concordances. In particular, we construct a large class of non-regular (and hence non-decomposable) Lagrangian concordances with Lagrangian fillable (and hence augmentable and non-stabilised) Legendrian ends. These are the first examples of non-regular (and hence non-decomposable) Lagrangian cobordisms between connected non-stabilised Legendrian submanifolds.  Our construction allows us to say that the Lagrangian concordance relation is not a partial order among Lagrangian fillable (and hence augmentable) Legendrian knots. 
More precisely, we prove the following.

\begin{figure}[h]
\includegraphics[height=4.3cm]{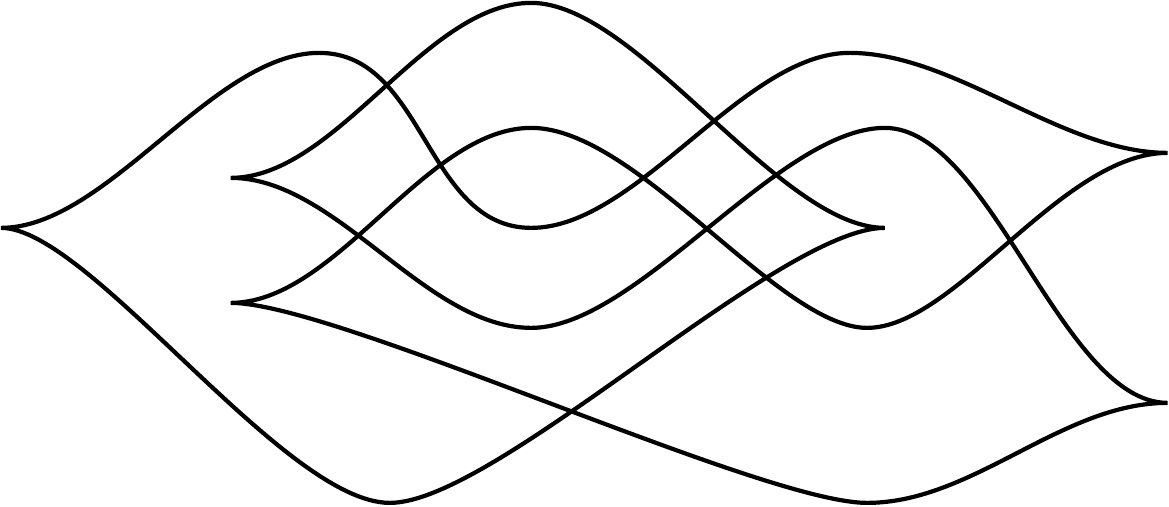}
\vspace{0.1in}
\caption{Front projection of $\Lambda_{m(9_{46})}$.}
\label{fig:m946}
\end{figure}

\begin{theorem}
\label{nonantisim}
Let $L$ be a decomposable Lagrangian concordance from $\Lambda_-$ to $\Lambda_+$, where $\Lambda_-$ and $\Lambda_+$ are Legendrian knots in two different smooth isotopy classes. 
There exists pairs of Legendrian knots $\Lambda'_-$, $\Lambda'_+$, where $\Lambda'_{\pm}$ is obtained from $\Lambda_{\pm}$ by a number $n \gg 0$ of double-stabilisations, followed by a sequence of Legendrian satellite operations, such that
\begin{enumerate}
  \item $\Lambda'_\pm$ are Lagrangian fillable and still in different smooth isotopy classes (e.g.~when the satellite construction is the Legendrian Whitehead double);
\item There is a pair of Lagrangian concordances $L^{+}_{-}$ from $\Lambda'_-$ to $\Lambda'_+$ and $L^{-}_{+}$ from $\Lambda'_+$ to $\Lambda'_-$; and
\item  $L^{-}_{+}$ is not  strongly homotopy-ribbon, and hence not regular, and also not decomposable. On the other hand, $L_{-}^+$ is regular.
\end{enumerate}
\end{theorem}
There are plenty of decomposable Lagrangian concordances between smoothly non-isotopic Legendrian knots to which our result can be applied. For instance, one can start with a decomposable Lagrangian slice disc of a knotted Legendrian, and remove a small unknotted disc in order to obtain a decomposable Lagrangian concordance from the unknot to a non-trivial knot.
\begin{ex}
The decomposable Lagrangian disc filling of the Legendrian representative $\Lambda_{m(9_{46})}$ of the mirror of $9_{46}$ knot from Figure \ref{fig:m946} that we considered in \cite{NonregStabilizedConc} gives rise to a decomposable concordance from the standard Legendrian unknot $U$ of ${\tt tb}=-1$ to $\Lambda_{m(9_{46})}$.
\end{ex}
We refer the reader to \cite{ObstrLagrConc, LegSatDecomLagrCob} for more examples.

Parts (1) and (2) of Theorem \ref{nonantisim} imply that the relation defined by Lagrangian concordances on Lagrangian fillable (and hence augmentable) Legendrian knots is not anti-symmetric, and therefore it leads to the following corollary.
\begin{corollary}
\label{cornotpartord}
The Lagrangian concordance relation does not define a partial order on the class of Lagrangian fillable (and hence augmentable) Legendrian knots.
\end{corollary}

\begin{remark}
The Legendrians $\Lambda'_\pm$ that we can construct always have  Lagrangian fillings of positive genus.  We do not know whether the method used in this paper allows to get Lagrangian fillings of genus zero.   Note that even though the Lagrangian concordance $L^-_+$ is not ribbon (since it is not even strongly homotopy-ribbon, see Section \ref{sec:ribbon_sthribbon})\color{black}, the Lagrangian filling of $\Lambda'_-$ obtained by concatenating the (positive genus) Lagrangian filling of $\Lambda'_+$ with the Lagrangian concordance is ribbon. This folows by a result due to Boileau--Orevkov \cite{BoileauOrevkovquasipositive}.
\end{remark}

The proof of Theorem \ref{nonantisim} relies on the variation of the recent result by Agol \cite{RibbonConcPartialOrdering} for strongly homotopy-ribbon concordances, which states that strongly homotopy-ribbon concordances form a partial order; see Theorem \ref{agolth} and Remark \ref{ribbonanti-symmetry} below. We apply this result to concordances that arise when applying the so-called Legendrian satellite operations; see Section \ref{sec:satellites_Whitehead doubles}  for the details. Note that  Breen recently proved that the regularity property for Lagrangian concordances is preserved under the action of Legendrian satellite operations; see \cite[Theorem 1.11]{BreenRegularSliceDecom}. When combined with the  variation of Agol's result (see Remark \ref{ribbonanti-symmetry} below), we conclude the following. Recall that Legendrian satellite operations preserve the Lagrangian concordance relation; namely, any Lagrangian concordance induces a Lagrangian ``satellite concordance'' between the corresponding Legendrian satellites. For a decomposable Lagrangian concordance with smoothly non-isotopic ends that admits a Lagrangian inverse concordance (which thus is non-regular by Agol's result), any Legendrian satellite operation that acts injectively on smooth knot types must preserve the non-regularity of the inverse concordance. This is the main mechanism used in the proof of Theorem \ref{nonantisim}.
\color{black}

\section{Legendrian satellites, stabilisations and Legendrian Whitehead doubles}
\label{sec:satellites_Whitehead doubles}

Now we discuss Legendrian satellites, stabilisations, Legendrian Whitehead doubles and their basic Legendrian and topological properties. 

The Legendrian satellite construction was first  defined by Ng in \cite{InvLegLinks} and by Ng--Traynor in \cite{NgSolTorusLinks}, \color{black}it can be seen as a generalization of the $n$-copy construction of Legendrian links studied by Michatchev in \cite{InvofLegandTransLinks}. Legendrian satellites has been a crucial tool that has been used both for providing new examples, as well as for constructing new invariants\color{black}; see for example \cite{EtnyreVertesiLegSat, LegSatDecomLagrCob, NgSolTorusLinks}.

\begin{figure}[!h]
\includegraphics[height=4.5cm]{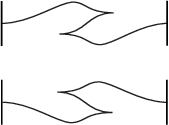}
\vspace{0.1in}
\caption{Front projection of $S_{+}$ (upper part) and of $S_-$ (lower part).}
\label{fig:stabilisations}
\end{figure}

\begin{figure}[!h]
\includegraphics[height=2cm]{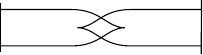}
\vspace{0.1in}
\caption{Front projection of $W_0$.}
\label{fig:tangle}
\end{figure}

Given a Legendrian knot $\Lambda\subset (S^3, \xi_{st})$ and a Legendrian $n$-tangle $\Pi\subset J^1([0,1])$, we define the Legendrian satellite of $\Lambda$ by $\Pi$ in the following way. First, we take  a canonical neighborhood $V(\Lambda)$ of $\Lambda$ contactomorphic to $J^1(S^1)$ and define the Legendrian satellite  of $\Lambda$ 
by $\Pi$ to be the link that is obtained from removing $V(\Lambda)$ and gluing in the solid torus containing $\Pi$. The Legendrian satellite of $\Lambda$ by $\Pi$ is typically denoted by $\Sigma(\Lambda, \Pi)$. When $\Pi$ equals $S_{\pm}$ from Figure~\ref{fig:stabilisations}, we say that $\Sigma(\Lambda, \Pi=S_{\pm})$ is a Legendrian $\pm$-stabilisation, the simplified notation of $\Sigma(\Lambda, \Pi=S_{\pm})$ is $S_{\pm}(\Lambda)$. \color{black} When $\Pi$  equals $W_0$ from Figure~\ref{fig:tangle}, we say that $\Sigma(\Lambda, \Pi=W_0)$ is a Legendrian Whitehead double. 

\subsection{Legendrian satellites and Lagrangian concordances}

Cornwell--Ng--Sivek in \cite[Theorem 2.4]{ObstrLagrConc} proved that Lagrangian concordance is preserved by all Legendrian satellite
operations. 
This is analogous to a well-known result in the classical knot concordance theory. 
In fact, after observing that Lagrangian cylinders admit standard neighborhoods, the argument is virtually identical. This result of Cornwell--Ng--Sivek will be important for us when $\Pi=W_0$ and $\Pi=S_{\pm}$.

\subsection{Properties of Legendrian Whitehead doubles}

\color{black}

\begin{figure}[!h]
\includegraphics[height=2cm]{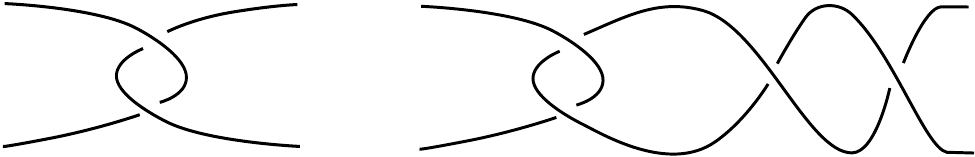}
\vspace{0.1in}
\caption{Smooth Whitehead doubling pattern $W_0^{C^{\infty}}$  (left) and smooth twisted  Whitehead doubling pattern $\Delta W_0^{C^{\infty}}$ (right).}
\label{fig:smoothWhitehead}
\end{figure}

We now make a few useful observations about Legendrian Whitehead doubles: 
\begin{itemize}
\item[(1)] Legendrian Whitehead doubles of Legendrian knots are obviously connected, i.e.  Legendrian Whitehead doubles of Legendrian knots are Legendrian knots.
\item[(2)] From the work of Bourgeois--Sabloff--Traynor \cite[Section 5.1, Proposition 5.1]{GeogrBourSabloffTraynor} it follows that 
for an arbitrary Legendrian knot $\Lambda$ in $(S^3, \xi_{st})$, its Legendrian
Whitehead double satisfies ${\tt rot}(\Sigma(\Lambda, W_0))=0$, ${\tt tb}(\Sigma(\Lambda, W_0))=1$ and admits a genus one Lagrangian filling. Besides that, if   
${\tt rot}(\Lambda)=0$, then $\Sigma(\Lambda, W_0)$
always admits a genus one gf-compatible Lagrangian filling. This, in particular, implies that $\Sigma(\Lambda, W_0)$ is augmentable and non-stabilised. \color{black}
\item[(3)] Finally, we recall that from the smooth point of view (twisted) Whitehead doubling can be seen as a well-defined operation on smooth isotopy classes of knots, which acts on  smooth isotopy classes of knots injectively. In other words, given two smoothly non-isotopic knots $K$ and $K'$, the smooth (twisted) Whitehead doubles of $K$ and $K'$ are also smoothly non-isotopic. We refer the reader to \cite{JSJdecomp} for this result. To that end, recall that the smooth isotopy class of a Legendrian Whitehead double  $\Sigma(\Lambda, W_0)$ does not always coincide with the smooth isotopy class of the corresponding smooth standard Whitehead double $\Sigma^{C^{\infty}}(\Lambda, W_0^{C^{\infty}})$. Rather, it coincides with the smooth isotopy class of the smooth \emph{twisted} Whitehead double $\Sigma^{C^{\infty}}(\Lambda, \Delta^{{\tt tb}(\Lambda)}W_0^{C^{\infty}})$, where $\Delta W_0^{C^{\infty}}$ denotes  the result of applying a full right handed twist to the pattern $W_0^{C^{\infty}}$; see Figure \ref{fig:smoothWhitehead}. For more details about this relation we refer the reader to \cite[Section 1]{EtnyreVertesiLegSat}.
\end{itemize}

\section{Inverting decomposable Lagrangian concordances}
Here we recall the construction from \cite{NonregStabilizedConc} that allows us to invert decomposable Lagrangian concordances. In the following we assume that we are given two Legendrian knots $\Lambda_{-}$, $\Lambda_+$  and a decomposable Lagrangian concordance $L$ from $\Lambda_-$ to $\Lambda_+$.

The construction is based on the following three ingredients:
\begin{itemize}
\item[(i)] Following \cite[Proposition 5.4]{NonregStabilizedConc} we consider the inverted cobordism
  $$L'=\{(-t,x,y,z)\in \R_x\times \R^3_{(x,y,z)}\ | \ (t,x,y,z)\in L\}.$$
 Even though the cobordism $L'$ is not necessarily Lagrangian, its tangent planes are homotopic to totally-real tangent planes, by the aforementioned proposition.  The $h$-principle for totally real embeddings thus implies that $L'$  admits a compactly supported smooth isotopy to a concordance that is totally real for any choice of cylindrical almost complex structure. Let us denote this totally real concordance from $\Lambda_+$ to $\Lambda_-$ by $T$.
\item[(ii)]  Then we apply  the approximation result of the first author from \cite{LagrApproxTotReal}. 
It guarantees that $T$ can be deformed to a Lagrangian
concordance $S(T)$ from $S(\Lambda_+)$ to $S(\Lambda_-)$ by a smooth isotopy, where both $S(\Lambda_\pm)$ are obtained from $\Lambda_{\pm}$
by adding a sufficiently large number of stabilisations  of both signs.
\item[(iii)] Finally, since $L$ is a Lagrangian concordance from $\Lambda_-$ to $\Lambda_+$, using either the result of Cornwell--Ng--Sivek \cite{ObstrLagrConc} or Chantraine \cite{LagrConcLegKnots} we get that there is a Lagrangian concordance $S(L)$ from $S(\Lambda_-)$ to $S(\Lambda_+)$.
\end{itemize}
All in all, these three ingredients prove the following statement, which is a slight generalization of \cite[Theorem 5.6]{NonregStabilizedConc}:
\begin{theorem}
\label{invertingconc}
Given two Legendrian knots $\Lambda_{-}$, $\Lambda_+$ and a decomposable Lagrangian concordance $L$ from $\Lambda_-$ to $\Lambda_+$. 
Then there are two Lagrangian concordances 
$S(L)$ from $S(\Lambda_-)$ to $S(\Lambda_+)$ and $S(T)$ from $S(\Lambda_+)$ to $S(\Lambda_-)$, where $S(\Lambda_\pm)$ are obtained from $\Lambda_\pm$ by sufficiently many stabilisations of both signs.
\end{theorem}

\section{Proof of Parts (1) and (2) of Theorem \ref{nonantisim}}
We take  a decomposable Lagrangian concordance $L$ from $\Lambda_-$ to $\Lambda_+$, where $\Lambda_-$ and $\Lambda_+$ are smoothly non-isotopic Legendrian knots. Then from Theorem \ref{invertingconc} it follows that there is a pair of Lagrangian concordances $S(L)$ from $S(\Lambda_-)$ to $S(\Lambda_+)$ and $S(T)$ from $S(\Lambda_+)$ to $S(\Lambda_-)$.  Since $\Lambda_-$ is not smoothly isotopic to $\Lambda_+$ and since $\Lambda_{\pm}$ is smoothly isotopic to $S(\Lambda_\pm)$, it follows that $S(\Lambda_-)$ is not smoothly isotopic to $S(\Lambda_+)$. Hence, we see that $S(\Lambda_-)$ is not Legendrian isotopic to $S(\Lambda_+)$.

In order to contradict the anti-symmetry of Lagrangian concordance in the class of fillable or augmentable Legendrian submanifolds we do the following. We perform a satellite operation
on $S(\Lambda_-)$ and $S(\Lambda_+)$ in such a way that both $\Sigma(S(\Lambda_{\pm}), \Pi)$ are fillable. This can be achieved in particular by taking Legendrian Whitehead doubles $\Sigma(S(\Lambda_{\pm}), W_0)$.  As we have seen in Observation (2) from Section \ref{sec:satellites_Whitehead doubles}, the work of Bourgeois--Sabloff--Traynor \cite[Proposition 5.1]{GeogrBourSabloffTraynor} implies that $\Sigma(S(\Lambda_{\pm}), W_0)$ both admit 
Lagrangian fillings of genus one, and hence, in particular, are augmentable and non-stabilised.  In addition, since $S(\Lambda_-)$ is not smoothly isotopic to $S(\Lambda_+)$, following the discussion in Observation (3) from Section \ref{sec:satellites_Whitehead doubles} we note that Legendrian Whitehead doubling acts injectively on smooth isotopy classes of knots, and hence  $\Sigma(S(\Lambda_{-}), W_0)$ is not smoothly isotopic (and thus not Legendrian isotopic) to $\Sigma(S(\Lambda_{+}), W_0)$. This leads to part (1) of Theorem~\ref{nonantisim}.

Now by applying the result of Cornwell--Ng--Sivek \cite{ObstrLagrConc} to $S(L)$ and $S(T)$ we get Lagrangian concordances  $\Sigma(S(L), W_0)$ from $\Sigma(S(\Lambda_{-}), W_0)$ to $\Sigma(S(\Lambda_{+}), W_0)$ and $\Sigma(S(T), W_0)$ from $\Sigma(S(\Lambda_{+}), W_0)$ to $\Sigma(S(\Lambda_{-}), W_0)$. This implies part (2) of Theorem~\ref{nonantisim}. 

Finally, as we have seen $\Sigma(S(\Lambda_{-}), W_0)$ is not Legendrian isotopic to $\Sigma(S(\Lambda_{+}), W_0)$, and hence the existence of Lagrangian concordances $\Sigma(S(L), W_0)$ and $\Sigma(S(T), W_0)$ contradicts anti-symmetry for the Lagrangian concordance relation on the class of fillable (or augmentable) Legendrian knots. Thus, we get Corollary~\ref{cornotpartord}.

\section{Stongly homotopy-ribbon concordances and the result of Agol}
\label{sec:ribbon_sthribbon}
Now we discuss the notions of ribbon and strongly homotopy-ribbon 
concordances of knots. These notions and the relations between them will be important for the proof of Part (3) of Theorem \ref{nonantisim}.

Given two knots $K_0$ and $K_1$ in $S^3$ and a concordance $C\subset [0,1]\times S^3$ from $K_0$ to $K_1$, we say that
\begin{itemize}
\item $C$ is {\em ribbon} if there is a compactly supported isotopy to a concordance whose projection to the $[0,1]$-factor is a Morse function with only critical points of index $0$ and $1$.
\item $C$ is {\em strongly homotopy-ribbon} if the complement of $C$ can be constructed by attaching $4$-dimensional $1$-handles
and $2$-handles to $S^3\setminus K_0$.
\end{itemize}

Following the discussion of Miller--Zemke in \cite[Section 1]{KFHstronglyhomrib}, 
we observe that  every ribbon concordance is strongly homotopy-ribbon. 
Note that  it is unknown whether there is a strongly homotopy-ribbon concordance that is not ribbon.

The following recent result of Agol, which was shown in \cite[Theorem 3.2]{ObstrLagrConc} by Cornwell--Ng--Sivek in the special case when one of the knots is the unknot, \color{black} is crucial for our construction:
 \begin{theorem}[\cite{RibbonConcPartialOrdering}]
\label{agolth}
Ribbon concordance relation of knots defines a partial order. 
\end{theorem}
\begin{remark}
\label{ribbonanti-symmetry}
Following the proof of \cite[Theorem 1.1]{RibbonConcPartialOrdering} and keeping in mind that every ribbon concordance is strongly homotopy-ribbon, we observe that the argument of \cite[Theorem 1.1]{RibbonConcPartialOrdering} holds for strongly homotopy-ribbon concordances. In other words, strongly homotopy-ribbon concordances define a partial order, and hence, in particular, strongly homotopy-ribbon concordances define an anti-symmetric relation. This remark also appears in the beginning of \cite[Section 3]{RibbonConcPartialOrdering}. A similar generalisation of \cite[Theorem 3.2]{ObstrLagrConc} to strongly homotopy-ribbon concordances also appeared in work \cite{NonregStabilizedConc} by the authors.
\end{remark}

\section{Proof of Part (3) of Theorem \ref{nonantisim}}

From Theorem \ref{nonantisim} (1)  we know that there is a pair of Lagrangian concordances $S(L)$ from $S(\Lambda_-)$ to $S(\Lambda_+)$ and $S(T)$ from $S(\Lambda_+)$ to $S(\Lambda_-)$, where $S(\Lambda_-)$ is not smoothly isotopic to $S(\Lambda_+)$.

Observe that since $L$ is decomposable, from the result of Conway--Etnyre--Tosun \cite[Theorem 1.10]{NonregImplNondecomp} it follows that $L$ is regular. 
Now from the result of Breen \cite[Theorem 1.11]{BreenRegularSliceDecom} it follows that $S(L)$ is regular. 
Here we use the fact  that taking a stabilisation can be seen as a particular case of taking a Legendrian satellite.

We claim that $S(T)$ is a non-regular Lagrangian concordance. To see this, we first observe that every regular Lagrangian concordance is strongly homotopy-ribbon. 
It comes from the following argument: if $L$ is a regular Lagrangian concordance in the symplectization of $S^3$ then, by the definition of regularity, we get that the ambient Weinstein structure can be homotoped so the Liouville vector field becomes tangent to $L$. This gives a Weinstein handle description of the complement of $L$. In dimension $4$, Weinstein handles have index $\le 2$, so the complement has a handle decomposition with no handles of index $>2$, which by definition means that $L$ is strongly homotopy-ribbon.

Now let us assume that $S(T)$ is a regular Lagrangian concordance. Then both $S(T)$ and $S(L)$ are regular and hence strongly homotopy-ribbon. Since  $\Lambda_-$ and $\Lambda_+$ are not smoothly isotopic,  $S(\Lambda_-)$ and $S(\Lambda_+)$ are smoothly non-isotopic. Then we get strongly homotopy-ribbon concordances $S(L)$, $S(T)$ from $S(\Lambda_-)$ to $S(\Lambda_+)$ and from $S(\Lambda_+)$ to $S(\Lambda_-)$. Since  $S(\Lambda_-)$ and $S(\Lambda_+)$ are smoothly non-isotopic, this leads to contradiction with the version of Agol's result \cite{RibbonConcPartialOrdering} for strongly homotopy-ribbon concordances; see Remark \ref{ribbonanti-symmetry}.

\begin{remark}
Recall that regular Lagrangian concordances and decomposable Lagrangian concordances define reflexive and transitive relations on smooth isotopy classes of knots in $S^3$.  
Since, as we discussed,  every regular Lagrangian concordance is strongly homotopy-ribbon, and since 
every decomposable concordance is ribbon,  Theorem \ref{agolth} and  Remark \ref{ribbonanti-symmetry} imply that 
\begin{itemize}
\item the relation defined by  regular Lagrangian concordances between Legendrian knots is a partial order on smooth isotopy classes of knots, and 
\item the relation defined by decomposable Lagrangian concordances  between Legendrian knots is a partial order on smooth isotopy classes of knots.
\end{itemize}
\end{remark}

Then, in order to extend this construction to fillable Legendrian submanifolds, we apply Legendrian Whitehead doubling. 
As we have already mentioned in the proof of Theorem \ref{nonantisim}, from the work of Bourgeois--Sabloff--Traynor \cite[Proposition 5.1]{GeogrBourSabloffTraynor} it follows that $\Sigma(S(\Lambda_{\pm}), W_0)$ both admit
Lagrangian fillings of genus one, and hence, in particular, are augmentable and not stabilised. Besides that, as we have already noted in Observation (3) from Section \ref{sec:satellites_Whitehead doubles}, since $S(\Lambda_{-})$ and $S(\Lambda_{+})$ are not smoothly isotopic, $\Sigma(S(\Lambda_{-}), W_0)$ and $\Sigma(S(\Lambda_{+}), W_0)$ are not smoothly isotopic.

Then we proceed with the same argument for Lagrangian concordances $\Sigma(S(L), W_0)$ and $\Sigma(S(T), W_0)$. More precisely, since $S(L)$ is regular, from the result of Breen \cite[Theorem 1.11]{BreenRegularSliceDecom} it follows that $\Sigma(S(L), W_0)$ is regular. We claim that $\Sigma(S(T), W_0)$  is a non-regular Lagrangian concordance. If we assume that  $\Sigma(S(T), W_0)$ is regular (and hence is strongly homotopy-ribbon), then we  get to  the contradiction, since we would have strongly homotopy-ribbon concordances $\Sigma(S(L), W_0)$, $\Sigma(S(T), W_0)$ between smoothly non-isotopic  $\Sigma(S(\Lambda_{-}), W_0)$, $\Sigma(S(\Lambda_{+}), W_0)$. This finishes the proof.

\begin{remark}
In the proof of Theorem  \ref{nonantisim} we can also rely on other Legendrian satellites $\Sigma(\cdot, \Pi)$ different from the Legendrian Whitehead double, as long as
\begin{itemize}
\item $\Sigma(\cdot, \Pi)$ distinguishes smooth isotopy classes of the Legendrian ends of the corresponding Lagrangian concordances. There are many examples of satellites of this type, we refer the reader to \cite{JSJdecomp}.
\item $\Sigma(\cdot, \Pi)$ applied to the negative Legendrian ends of the corresponding Lagrangian concordances are fillable.
\end{itemize}
An example of a different Legendrian satellite that satisfies these properties is an iterated (or multiple) Legendrian Whitehead double.
\end{remark}

\begin{figure}[t]
\includegraphics[height=4cm]{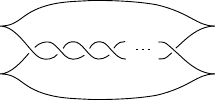}
\vspace{0.1in}
\caption{Front projection of the Legendrian torus $(2,n)$-knot discussed in \cite{LegKnotsLagrCob}. The total number of crossings is $n$.}
\label{fig:Lagr_T(2,n)}
\end{figure}

\begin{remark}
One can also produce knots $\Lambda''_\pm$ which are smoothly non-isotopic, Lagrangian fillable for some genus-$g$ surface for arbitrary choices of $g\ge 1$, and for which there exist Lagrangian concordances going both ways. One of these concordances will moreover be regular by construction, and thus the other one is necessarily non-regular by Agol's result. We can start by performing a Legendrian cusp-connected sum with the Legendrians $\Lambda'_\pm$ produced in the proof of Theorem \ref{nonantisim} and a Legendrian knot $\tilde{\Lambda}$ which admits a genus $g-1$ Lagrangian filling. For example, we can take $\tilde{\Lambda}$ to be one of the Legendrian torus $(2,n)$-knots from the work of Ekholm--Honda--Kalman \cite[Theorem 1.6]{LegKnotsLagrCob} shown in Figure \ref{fig:Lagr_T(2,n)}. Since $\Lambda'_\pm$ admits a Lagrangian filling of genus one, $\tilde{\Lambda}$ admits a Lagrangian filling of genus $g-1$ and since there exists a Lagrangian pair-of-pants cobordism from $\Lambda'_\pm\sqcup \tilde{\Lambda}$ to $\Lambda'_\pm\# \tilde{\Lambda}$, it follows that $\Lambda'_\pm\# \tilde{\Lambda}$ admits a Lagrangian filling of genus $g$. For the existence of the Lagrangian concordances from $\Lambda'_-\# \tilde{\Lambda}$ to $\Lambda'_+\# \tilde{\Lambda}$ and from $\Lambda'_+\# \tilde{\Lambda}$ to $\Lambda'_-\# \tilde{\Lambda}$, note that the Lagrangian concordance relation is preserved under the  action of cusp-connected sum; indeed, one can take a connected sum of the concordances $L_-^+$, $L_+^-$ from Theorem \ref{nonantisim} and the trivial concordance $[0,1]\times \tilde{\Lambda}$. Since $L_-^+$ and $[0,1]\times \tilde{\Lambda}$ are decomposable, their connected sum $L_-^+ \# ([0,1]\times \tilde{\Lambda})$ is decomposable, and hence by the result of Conway--Etnyre--Tosun  \cite[Theorem 1.10]{NonregImplNondecomp} the concordance is also regular. In addition, we note that since  $\Lambda'_\pm$ are not smoothly isotopic, from Schubert's theorem \cite{BZKnots} it follows that $\Lambda'_\pm\# \tilde{\Lambda}$ are not smoothly isotopic as well. Finally, using that $L_-^+ \# ([0,1]\times \tilde{\Lambda})$ is strongly homotopy-ribbon and applying the result of Agol, see Theorem \ref{agolth} and Remark \ref{ribbonanti-symmetry}, we get that $L_+^- \# ([0,1]\times \tilde{\Lambda})$ is non-regular, and hence is non-decomposable.
\end{remark}

\end{document}